\documentclass[11pt]{article}

\usepackage{float}
\usepackage{amsmath}
\usepackage{graphicx}
\usepackage{subfigure}
\usepackage[numbers]{natbib}
\usepackage{float}
\usepackage{geometry}
\usepackage{enumerate}
\usepackage{tikz}
\usepackage[english]{babel}
\usepackage{longtable}
\usepackage{color}
\usepackage{amssymb, amsmath, amsthm}
\usepackage{multirow}
\usepackage{paralist}
\usepackage[titletoc, title]{appendix}
\usepackage{authblk}
\usepackage{setspace}
\usepackage{dsfont}
\usepackage{enumerate}
\RequirePackage[OT1]{fontenc}
\usepackage{amsthm,amsmath}

\RequirePackage[colorlinks,citecolor=blue,urlcolor=blue]{hyperref}

\newtheorem{theorem}{Theorem}
\newtheorem{lemma}{Lemma}

\DeclareMathOperator{\expit}{expit}
\DeclareMathOperator{\logit}{logit}

\newcommand{\dd}{\mathrm{d}}
\newcommand{\one}{\mathds{1}}

\newcommand{\R}{\mathbb{R}}

\DeclareMathOperator*{\argmin}{\arg\!\min}
\title{Second-Order Inference for the Mean of a Variable Missing at Random}
\author[1]{Iv\'an  D\'iaz\thanks{corresponding author: idiaz@google.com}}
\author[2]{Marco Carone}
\author[3]{Mark J. van der Laan}
\affil[1]{\small Google Inc.}
\affil[2]{\small Department of Biostatistics, University of Washington.}
\affil[3]{\small Division of Biostatistics, University of California
  at Berkeley.}

\begin{document}\maketitle

\begin{abstract}
  We present a second-order estimator of the mean of a variable subject to
  missingness, under the missing at random assumption. The estimator improves
  upon existing methods by using an approximate second-order expansion of the
  parameter functional, in addition to the first-order expansion employed by
  standard doubly robust methods. This results in weaker assumptions about the
  convergence rates necessary to establish consistency, local efficiency, and
  asymptotic linearity. The general estimation strategy is developed under the
  targeted minimum loss-based estimation (TMLE) framework. We present a
  simulation comparing the sensitivity of the first and second order estimators
  to the convergence rate of the initial estimators of the outcome regression
  and missingness score. In our simulation, the second-order TMLE improved the
  coverage probability of a confidence interval by up to $85\%$. In addition, we
  present a first-order estimator inspired by a second-order expansion of the
  parameter functional. This estimator only requires one-dimensional smoothing,
  whereas implementation of the second-order TMLE generally requires kernel
  smoothing on the covariate space. The first-order estimator proposed is
  expected to have improved finite sample performance compared to existing
  first-order estimators. In our simulations, the proposed first-order estimator
  improved the coverage probability by up to $90\%$. We provide an illustration
  of our methods using a publicly available dataset to determine the effect of
  an anticoagulant on health outcomes of patients undergoing percutaneous
  coronary intervention. We provide R code implementing the proposed estimator.
\end{abstract}

\section{Introduction}\label{intro}

Estimation of the mean of an outcome subject to missingness has been extensively
studied in the literature. Under the assumption that missingness is independent
of the outcome conditional on observed covariates, the marginal expectation is
identified as a parameter depending on the conditional expectation given
covariates among observed individuals (outcome regression henceforth) and the
marginal distribution of the covariates. If the covariate vector consists of a
few categorical variables, a nonparametric maximum likelihood estimator yields
an optimal (i.e., asymptotically efficient) estimator of the mean
outcome. However, if the covariate vector contains continuous variables or its
dimension is large, estimation of the outcome regression requires smoothing on
the covariate space. This has often been achieved by means of a parametric
model. Unfortunately, the correct specification of a parametric model is a
chimerical task in high-dimensional settings or in the presence of continuous
variables \cite{starmans2011models}, and data-adaptive estimation methods such
as those developed in the statistical learning literature (e.g., super learning,
model stacking, bagging) must be used.

Our methods are developed in the context of targeted learning
\cite{vanderLaanRose11,vanderLaan&Rubin06}, a branch of statistics that deals
with the use of data-adaptive methods coupled with optimal estimation theory for
infinite-dimensional models. In particular, the targeted minimum loss-based
estimation (TMLE) framework allows consistent and locally efficient estimation
of arbitrary low-dimensional parameters in high-dimensional models under
regularity and smoothness conditions. In our context, targeted learning allows
the incorporation of flexible data-adaptive estimators of the outcome regression
into the estimation procedure.

Several doubly robust and locally efficient estimators have been proposed for
the missing data problem. These estimators are based on a first-order expansion
of the parameter functional, and are asymptotically efficient, under certain
conditions. Arguably, the most important condition is that the outcome
regression and the probability of missingness conditional on covariates
(missingness score henceforth) are estimated consistently at an appropriate
rate. A sufficient assumption for establishing $\sqrt{n}$-consistency of doubly
robust estimators is that the outcome regression and the missingness score
converge to their true values at rates faster than $n^{-1/4}$. In this paper we
are concerned with asymptotically efficient estimation under slower consistency
rates of these estimators. In particular, we present a second-order
TMLE that incorporates a second-order expansion of the parameter functional in
order to relax this assumption, which may be implausible in high dimensions and
for certain data-adaptive estimators. The method we present is an application of
the general higher-order estimation theory we present in \cite{carone2014}. We
refer to the second-order estimator as 2-TMLE in contrast to the first-order
TMLE discussed by \cite{vanderLaan&Rubin06}, referred to as 1-TMLE.

A complete literature review of higher-order estimation theory is presented in
\cite{carone2014}. The most relevant references for the problem studied here are
\cite{robins2009quadratic} and \cite{robins2009semiparametric}. In particular,
\cite{robins2009quadratic} presents a particular second-order expansion of the
target parameter, as well as a second-order estimator based on that
expansion. This estimator directly uses inverse weighting by a kernel estimate
of the covariate density. As a result of the curse of dimensionality, the
estimator may perform poorly in finite samples as the dimension of the covariate
space increases. Particularly, it may fall outside of the parameter space. In
contrast, the 2-TMLE presented here is a substitution estimator that always
falls in the parameter space. The results presented in \cite{tan2014second}
establish the asymptotic properties of various calibration estimators in the
context of missing data problems, concluding that some of them are second-order
estimators. However, their results are not directly related to this manuscript
since they assume a Euclidean parametrization of the outcome model and a known
missingness score. 

As with the estimator presented in \cite{robins2009quadratic}, implementation of
the 2-TMLE requires approximating the second-order influence function by means
of kernel smoothing. When the covariate space is high-dimensional,
this approximation is subject to the curse of dimensionality. This issue
may be circumvented by utilizing an alternative second-order expansion that uses
kernel smoothing on the missingness score, which is a one-dimensional function
of the covariate vector. Since the true missingness score is generally unknown,
implementation of this estimator must be carried out using an estimated
missingness score. Unfortunately, introduction of the estimated missingness
score in place of its true value yields a second-order remainder term in the
analysis of the estimator. As a consequence, the estimator obtained is not a
second-order estimator. We refer to this estimator as a $1^\star$-TMLE in
accordance with this observation. Notably, the second-order remainder term
obtained with the $1^\star$-TMLE is different from that of the 1-TMLE, which
implies they have different finite sample properties. We conjecture that the
$1^\star$-TMLE improves finite sample performance over the 1-TMLE, and present a
case study in which there are considerable finite sample gains.

Compared to the standard 1-TMLE, implementation of the $1^\star$-TMLE requires
the inclusion of one additional covariate in the outcome regression. As a
result, its implementation is straightforward and comes at no computational
cost. Moreover, the potential finite sample gains in performance can be
overwhelming, as we illustrate in a simulation studying the coverage probability
and mean squared error of the two estimators.

The paper is organized as follows. In Section~\ref{foet} we review first-order
efficient estimation theory for the mean outcome in a missing data model. In
Section~\ref{soet} we present the second-order expansion of the parameter
functional and use it in Section~\ref{estimaw} to construct a 2-TMLE. In
Section~\ref{estimag} we introduce the $1^\star$-TMLE discussed
above. Section~\ref{simula} presents a simulation showing that the
$1^\star$-TMLE and the 2-TMLE have improved coverage probabilities and mean
squared error for slow convergence rates of the estimated outcome regression and
missingness score. We conclude with Section~\ref{aplica} illustrating the use of
the $1^\star$-TMLE in a real data application.

\section{Review of First-Order Estimation Theory}\label{foet}
Let $W$ denote a $d$-dimensional vector of covariates, and let $Y$ denote an
outcome of interest measured only when a missingness indicator $A$ is equal to
one. To simplify the exposition, we assume that $Y$ is binary or continuous
taking values in the interval $(0,1)$. The observed data $O=(W,A,AY)$ is assumed
to have a distribution $P_0$ in the nonparametric model $\cal M$. Assume we
observe an i.i.d. sample $O_1,\ldots,O_n$, and denote the empirical distribution
by $P_n$. For every element $P\in \cal M$, we define
\begin{align*}
  Q_W(P)(w)&:= P(W\leq w)\\
  g(P)(w)&:= P(A=1|W=w)\\
  \bar Q(P)(w) &:= E_P(Y|A=1, W=w),
\end{align*}
where $E_P$ denotes expectation under $P$. We denote $Q_{W,0}:=Q_W(P_0)$,
$g_0:=g(P_0)$, and $\bar Q_0:=\bar Q(P_0)$. We refer to $\bar Q$ as the
\textit{outcome regression}, and to $g$ as the \textit{missingness score}. We
suppress the argument $P$ from the notation $Q_W(P)$, $g(P)$, and $\bar Q(P)$
whenever it does not cause confusion. For a function $f$ of $o$, we use the
notation $Pf:=\int f(o)\dd P(o)$. Let $\Psi:\cal M\to \R$ be a parameter mapping
defined as $\Psi(P):=E_P\{\bar Q(W)\}$, and let $\psi_0:=\Psi(P_0)$. Under the
assumptions that missingness $A$ is independent of the outcome $Y$ conditional
on the covariates $W$ and $P_0(g_0(W)>0)=1$, it can be shown that
$\psi_0=E_{F_0}(Y)$, where $F_0$ is the true distribution of the full data
$(W,Y)$. Because $\Psi$ depends on $P$ only through $Q:=(Q_W, \bar Q)$, we also
use the alternative notation $\Psi(Q)$ to refer to $\Psi(P)$.

First-order inference for $\psi_0$ is based on the following expansion of the
parameter functional $\Psi(P)$ around the true $P_0$:
\begin{equation}
  \Psi(P) - \Psi(P_0) = -P_0 D^{(1)}(P) + R_2(P,P_0),\label{foexp}
\end{equation}
where $D^{(1)}(P)$ is a function of an observation $o=(w,a,y)$ that depends on
$P$, and $R_2(P,P_0)$ is a second-order remainder term. The super index $(1)$ is
used to denote a first-order approximation. This expansion may be seen as
analogous to a Taylor expansion when $P$ is indexed by a finite-dimensional
quantity, and the expression \textit{second-order} may be interpreted in the
same way.

We use the expression \textit{first-order estimator} to refer to estimators
based on first-order approximations as in equation
(\ref{foexp}). Analogously, the expression \textit{second-order
  estimator} is used to refer to estimators based on second-order
approximations, e.g., as presented in Section \ref{soet} below.

Doubly robust locally efficient inference is based on approximation (\ref{foexp}) with
\begin{align}
  D^{(1)}(P)(o)&=\frac{a}{g(w)}\{y - \bar Q(w)\} + \bar Q(w) -
  \Psi(P),\label{defD1}\\
  R_2(P,P_0) &= \int \left\{1-\frac{g_0(w)}{g(w)}\right\}\{\bar Q(w) -
  \bar Q_0(w)\}\dd Q_{W,0}(w).\label{defR2}
\end{align}
Straightforward algebra suffices to check that equation
(\ref{foexp}) holds with the definitions given above. $D^{(1)}$ as defined in
(\ref{defD1}) is referred to as the canonical gradient or the  efficient
influence function \cite{bang2005doubly,vanderLaanRose11}.

First-order targeted minimum loss-based estimation
of $\psi_0$ is performed in the following steps:
\begin{enumerate}[{Step }1.]
\item \textit{Initial estimators.} Obtain initial estimators  $\hat g$
  and $\hat{\bar Q}$ of  $g_0$ and $\bar Q_0$. In general, the
  functional form of $g_0$ and $\bar Q_0$ will be unknown to the researcher.
  Since consistent estimation of these quantities is key to achieve asymptotic
  efficiency  of $\hat \psi$, we advocate for the use of data-adaptive
  predictive  methods that allow flexibility in the specification of these
  functional forms.
\item \textit{Compute auxiliary covariate.} For each
  subject $i$, compute the auxiliary covariate \[\hat H^{(1)}(W_i) :=
  \frac{1}{\hat g(W_i)}.\]
\item \textit{Solve estimating equations.} Estimate the parameter
  $\epsilon$ in the logistic regression model
  \begin{equation}
    \logit \hat{\bar Q}_{\epsilon,h}(w) = \logit \hat{\bar Q}(w) + \epsilon
    \hat H^{(1)}(w),\label{submodel1}
  \end{equation}
  by fitting a standard logistic regression model of $Y_i$ on $\hat
  H^{(1)}(W_i)$, with no intercept and with
  offset $\logit \hat{\bar Q}(W_i)$, among observations with
  $A=1$. Alternatively, fit the model
  \begin{equation*}
    \logit \hat{\bar Q}_{\epsilon,h}(w) = \logit \hat{\bar Q}(w) + \epsilon
  \end{equation*}
  with weights $\hat H^{(1)}(W_i)$ among observations with
  $A=1$. In either case, denote the estimate of $\epsilon$ by $\hat \epsilon$.
\item \textit{Update initial estimator and compute 1-TMLE.} Update the
  initial estimator as $\hat{\bar Q}^\star_h(w) = \hat{\bar Q}_{\hat
    \epsilon}(w)$, and define the 1-TMLE as $\hat\psi=\Psi(\hat{\bar Q}^\star)$.
\end{enumerate}
Note that this estimator $\hat P$ of $P_0$ satisfies $P_n D^{(1)}(\hat
P)=0$. For a full presentation of the TMLE algorithm the interested reader is
referred to \cite{vanderLaanRose11} and the references therein. Using equation
(\ref{foexp}) along with  $P_n D^{(1)}(\hat P)=0$ we obtain that
\[\hat \psi-\psi_0=(P_n-P_0)D^{(1)}(\hat P) + R_2(\hat P, P_0).\]
Provided that
\begin{enumerate}[i)]
\item $D^{(1)}(\hat P)$ converges to $D^{(1)}(P_0)$ in $L_2(P_0)$ norm, and
\item the size of the class of functions considered for estimation of
  $\hat P$ is bounded (technically, there exists a Donsker class $\cal H$ so that
  $D^{(1)}(\hat P)\in \cal H$ with probability tending to one),
\end{enumerate}
results from empirical process theory (e.g., theorem 19.24 of
\cite{vanderVaart98}) allow us to conclude that
\[\hat\psi-\psi_0=(P_n-P_0)D^{(1)}(P_0) + R_2(\hat P, P_0).\]
In addition, if
\begin{equation}
  R_2(\hat P, P_0)=o_P(1/\sqrt{n}),\label{Rrootn}
\end{equation}
we obtain that $\hat\psi-\psi_0=(P_n-P_0)D^{(1)}(P_0) + o_P(1/\sqrt{n})$. This
implies, in particular, that $\hat\psi$ is a $\sqrt{n}$-consistent estimator of
$\psi_0$, it is asymptotically normal, and it is locally efficient.

In this paper we discuss ways of constructing an estimator that requires a
consistency assumption weaker than (\ref{Rrootn}). Note that (\ref{Rrootn}) is
an assumption about the convergence rate of a second order term involving the
product of the differences $\hat{\bar Q}-Q_0$ and $\hat g-g_0$. Using the
Cauchy-Schwarz inequality repeatedly, $|R_2(\hat P, P_0)|$ may be bounded as
\[|R_2(\hat P, P_0)|\leq ||1/\hat g||_{\infty}\, ||\hat g -
g_0||_{P_0}\,||\hat{\bar Q} - \bar Q_0||_{P_0},\] where $||f||^2_{P} := \int
f^2(o)\dd P(o)$, and $||f||_{\infty}:=\sup\{f(o):o\in \cal O\}$. For assumption
(\ref{Rrootn}) to hold is, it is sufficient to have that
\begin{enumerate}[i)]
\item $\hat g$ is bounded away from zero with probability tending to
  one;
\item $\hat g$ is the MLE of $g_0\in {\cal G}= \{ g(w;\beta):\beta \in
  \R^d\}$ (i.e., $g_0$ is estimated in a correctly specified
  parametric model) since this implies $||\hat g - g_0||_{P_0}=O_P(1/\sqrt{n})$;
  and
\item $||\hat{\bar Q} - \bar Q_0||_{P_0}=o_P(1)$.
\end{enumerate}
Alternatively the roles of $\hat g$ and $\hat{\bar Q}$ could also be
interchanged in i) and ii). As discussed in \cite{starmans2011models}, however,
correct specification of a parametric model is hardly achievable in
high-dimensional settings. Data-adaptive estimators must then be used for the
outcome regression and missingness score, but they may potentially yield a
remainder term $R_2$ with a convergence rate slower than $n^{-1/2}$. In the next
section we present a second-order expansion of the parameter functional that
allows the construction of estimators that require consistency assumptions
weaker than (\ref{Rrootn}).

\section{Second-Order Estimation}\label{soet}

Let us first introduce some notation. For a function $f^{(2)}$ of a pair of
observations $(o_1,o_2)$, let $P_0^2f^{(2)}:=\iint
f^{(2)}(o_1,o_2)\dd P_0(o_1)\dd P_0(o_2)$ denote the expectation of
$f^{(2)}$ with respect to the product measure $P_0^2$.

Second-order estimators are based on second-order expansions of the
parameter functional of the form
\begin{equation}
  \Psi(P)-\Psi(P_0)=-P_0D^{(1)}(P) - \frac{1}{2}P_0^2 D^{(2)}(P) +
  R_3(P,P_0),\label{soexp1}
\end{equation}
where $D^{(2)}(P)$ is a function of a pair of observations $(o_1,o_2)$ that
depends on $P$, and $R_3(P,P_0)$ is a third-order remainder term. $D^{(2)}$ is
referred to as a second-order gradient. This representation exists only if $W$
has finite support. If the support of $W$ is infinite, it is necessary to use an
approximate second-order influence function relying on smoothing, which yields a
bias term referred to as the \textit{representation error}. This may introduce
challenges due to the curse of dimensionality. In this section we discuss two
possible estimation strategies: (i) an estimator that implements kernel
smoothing on the covariate vector, and (ii) an estimator that implements kernel
smoothing on the missingness score. Strategy (i) is only practical in the
presence of a few, possibly data-adaptively selected covariates, although a
greater number of covariates may be included as sample size increases. Strategy
(ii) requires a-priori knowledge of the true missingness score, and is therefore
not applicable in most practical situations. As a solution, we propose to use
strategy (ii) with the estimated missingness score to obtain an estimator we
refer to as $1^\star$-TMLE. As discussed below, the $1^\star$-TMLE is not a
second-order estimator, since introduction of an estimated missingness score
yields a second-order term in the remainder term. Nevertheless, the potential
finite sample gains obtained with the $1^\star$-TMLE compared to the standard
1-TMLE are worth further investigation. In Section~\ref{simg} we present a
simulation study in which the $1^\star$-TMLE showed considerable finite sample
improvement in both mean squared error and coverage probability of associated
confidence intervals.

\subsection{Second-Order Estimator with Kernel Smoothing on the
  Covariate Vector}\label{estimaw}

Assume $W$ contains only discrete variables. Then the second
order expansion (\ref{soexp1}) holds with
\begin{align*}
  D^{(2)}(P)(o_1,o_2) &=
  \frac{2a_1\one\{w_1=w_2\}}{g(w_1)q_W(w_1)}\left\{1-\frac{a_2}{g(w_1)}\right\}\{y_1
  - \bar Q(w_1)\},\\ R_3(P,P_0)&=\int\left\{1-\frac{g_0(w)q_{W,0}(w)}{g(w)q_W(w)}\right\}\left\{1-\frac{g_0(w)}{g(w)}\right\}\left\{\bar
    Q(w)-\bar Q_0(w)\right\}\dd Q_{W,0}(w),
\end{align*}
where $q_W$ denotes the probability mass function associated to $Q_W$, and
$D^{(1)}$ is defined in (\ref{defD1}). It is easy to explicitly check that
equation (\ref{soexp1})  holds.

In most practical situations, however, $W$ is high-dimensional or it
contains continuous variables so that the indicator $\one\{w_1=w_2\}$
is essentially always zero. To circumvent this issue, we propose to use the
above expansion replacing the indicator function with a kernel function
$K_h(w_1-w_2)$ for a given bandwidth $h$. If $W$ takes values on a discrete
set, we define $K_h(w)=\one(w=0)$, so that the estimator $\hat g_h$ below is the
non-parametric estimator using empirical means in strata defined by $W$.
We denote the corresponding approximation of $D^{(2)}$ by $D_h^{(2)}$. The
following lemma establishes conditions under which
the representation error is negligible.

\begin{lemma}\label{lemmasmooth}
Suppose that the distribution of $W$ has compact support and is absolutely
continuous with respect to Lebesgue measure with density $Q_{W,0}$. Suppose that
$\hat Q_W$ is a working estimate of $Q_{W,0}$. If \begin{enumerate}
\item both $g_0$ and $Q_{W,0}$ are $(m_0+1)$-times continuously differentiable
  almost surely;\vspace{-.1in}
\item $K$ is orthogonal to all polynomial powers up until $m_0$;\vspace{-.1in}
\item there exists some $\delta>0$ such that $g_0$ is bounded below by
  $\delta$, and both $\hat g$ and $\hat Q_W$ are bounded below by $\delta$  with
  probability tending to one,
\end{enumerate} then we have that
\[
P_0^2 D^{(2)}_h(\hat{\bar Q}^\star,\hat g, \hat Q_W)-\lim_{h\rightarrow
  0}P_0^2 D^{(2)}_h(\hat{\bar Q}^\star,\hat g,\hat
Q_W)=O_P\left(h^{m_0+1}\|\hat{\bar Q}^\star-\bar{Q}_0\| \right),\]
where $\|\hat{\bar Q}^\star-\bar{Q}_0\|^2:=\int (\hat{\bar
  Q}^\star-\bar{Q}_0)^2(w)\dd Q_{W,0}(w)$.
\end{lemma}
The result above explicitly deals with kernel smoothing with common
bandwidth in all dimensions. The lemma also holds, however, if a
multivariate bandwidth is utilized, with $h$ substituted by $\max_jh_j$ in the
statement of the lemma.

\subsubsection{A corresponding 2-TMLE}
Analogous to the 1-TMLE discussed in the previous section, we construct
an estimator $\hat P$ satisfying $P_n D^{(1)}(\hat P) = P_n^2
D_h^{(2)}(\hat P) = 0$. Solving these equations allows us to exploit expansion
(\ref{soexp1}) and construct a $\sqrt{n}$-consistent estimator in
which assumption $R_2(\hat P,P_0)=o_P(1/\sqrt{n})$ is replaced by the
weaker assumption $R_3(\hat P,P_0)=o_P(1/\sqrt{n})$.

For a fixed bandwidth $h$, the proposed 2-TMLE is given by the following
algorithm, which is implemented in the R code provided in the supplementary
material.

\begin{enumerate}[{Step }1.]
\item \textit{Initial estimators.} See the previous section on the 1-TMLE.
\item \textit{Compute auxiliary covariates.} For each
  subject $i$, compute auxiliary covariates
  \begin{align*}
    \hat H^{(1)}(W_i)&:= \frac{1}{\hat g(W_i)}\\
    \hat H^{(2)}_h(W_i)&:= \frac{1}{\hat g(W_i)}\left\{1 -
      \frac{\hat g_h(W_i)}{\hat g(W_i)}\right\},\\
  \end{align*}
  where
  \[\hat g_h(w)=\frac{\sum_{i=1}^nK_h(w - W_i)A_i}{\sum_{i=1}^nK_h(w -
    W_i)}\]
  is a kernel regression estimator of $g_0(w)$.
\item \textit{Solve estimating equations.} Estimate the parameter
  $\epsilon=(\epsilon_1,\epsilon_2)$ in the logistic regression model
  \begin{equation}
    \logit \hat{\bar Q}_{\epsilon,h}(w) = \logit \hat{\bar Q}(w) + \epsilon_1
    \hat H^{(1)}(w) + \epsilon_2 \hat H^{(2)}_{h}(w),\label{submodel}
  \end{equation}
  by fitting a standard logistic regression model of $Y_i$ on $\hat
  H^{(1)}(W_i)$ and $\hat{\bar H}^{(2)}_{h}(W_i)$, with no intercept and with
  offset $\logit \hat{\bar Q}(W_i)$, among observations with $A=1$. Denote the
  estimate of $\epsilon$ by $\hat \epsilon$.
\item \textit{Update initial estimator and compute 2-TMLE.} Update the
  initial estimator as $\hat{\bar Q}^\star_h(w) = \hat{\bar Q}_{\hat
    \epsilon, h}(w)$, and define the $h$-specific 2-TMLE as
  $\hat\psi_h=\Psi(\hat{\bar Q}^\star_h)$
\end{enumerate}

The estimators presented above required a user-selected bandwidth $h$. Here we
discuss briefly two possible ways to select a bandwidth $\hat h$ to use in
practice. Certain convergence rates are required of this bandwidth so that the
resulting estimators achieve second-order properties (see Theorem
\ref{thsecondordertmleex} below). The first and easiest option is to select the
bandwidth that maximizes the log-likelihood loss function of the density
$q_0$. However, because this choice is targeted to estimation of $q_0$, it may
be sub-optimal for estimation of $\psi_0$. The second alternative is to use
the collaborative TMLE (C-TMLE) presented in \cite{vanderLaan&Gruber09}, which
may result in correct convergence rates as argued in \cite{carone2014}. The
question of whether these selectors achieve the required convergence rate is an
open research problem and will be the subject of future research.

The theorem below provides the exact conditions that guarantee asymptotic
linearity of $\hat \psi$.

\begin{theorem}\label{thsecondordertmleex}

Under the conditions of Lemma \ref{lemmasmooth}, and provided that
\begin{enumerate}
\item each of $\hat g-g_0$, $\hat{\bar Q}^\star-\bar{Q}_0$ and $\hat Q_W$ tend
 to zero in $L^2(Q_{W,0})$-norm;\vspace{-.1in}
 \item  there exists some $\delta>0$ such that $g_0$, $\hat g$ and $\hat
   Q_W\cdot \hat g$ are bounded below by $\delta$ with probability tending to
   one;\vspace{-.1in}
 \item each of $\hat g$, $\hat{\bar Q}^\star$ and $\hat Q_W$ have uniform
   sectional variation norm bounded by some $M<\infty$ with probability tending
   to one;\vspace{-.1in}
\item the kernel function $K$ is $2d$-times differentiable and $\hat
  h^{2d}n\rightarrow +\infty$,
\end{enumerate}
and either of
\begin{enumerate}\setcounter{enumi}{4}
\item[5a.] $R_2(\hat P,P_0)=o_P(n^{-1/2})$; or,\vspace{-.1in}
\item[5b.]\ $R_3(\hat P,P_0)=o_P(n^{-1/2})$\ \ \mbox{and}\ \ $\|\hat{\bar
    Q}^\star-\bar{Q}_0\| \hat h^{m_0+1}=o_P(n^{-1/2})$
\end{enumerate}
holds, $\hat \psi_{\hat h}$ is an asymptotically efficient estimator of $\psi_0$.
\end{theorem}

The proof of this theorem is presented in the supplementary materials.
A key argument in the proof is that $\hat P$ solves the estimating equations
$P_nD^{(1)}(\hat P)=P_n^2D^{(2)}_h(\hat P)=0$. The score
equations of the logistic regression model (\ref{submodel}) are equal to
\[  \sum_{i=1}^n\hat H^{(1)}(Y_i - \hat{\bar Q}_{\epsilon,h}(W_i))=0\quad\text{ and }\quad
  \sum_{i=1}^n\hat H^{(2)}_h(Y_i - \hat{\bar Q}_{\epsilon,h}(W_i))=0.\]
Because the maximum likelihood estimator solves the score equations,
it can be readily seen that
\[  \sum_{i=1}^n\hat H^{(1)}(Y_i - \hat{\bar Q}_h^\star(W_i))=0\quad\text{ and }\quad
  \sum_{i=1}^n\hat H^{(2)}_h(Y_i - \hat{\bar Q}_h^\star(W_i))=0,
\]
which, from the definitions of $\hat H^{(1)}$ and $\hat H^{(2)}_h$, correspond to
$P_nD^{(1)}(\hat P)=0$ and $P_n^2D^{(2)}_h(\hat P)=0$,
respectively.

As is evident from the conditions of the theorem, the rate at which the
bandwidth $\hat h$ decreases plays a critical role in the asymptotic behavior of
the $2$-TMLE described. On one hand, condition 5b of the theorem requires that
the bandwidth converge to zero sufficiently quickly in order for $n^{1/2}\|
\hat{\bar Q}^\star-\bar{Q}_0\| h^{m_0+1}$ to itself converge to zero, where
$m_0$ is the order of the kernel $K$ used. This ensures that the representation
error is negligible. On the other hand, condition 4 requires $\hat h$ to
converge to zero slowly enough to allow control of a $V$ statistic term
displayed in the proof of the theorem in the appendix.


Scrutiny of the theorem above reveals that a 2-TMLE will indeed generally be
asymptotically linear and efficient in a larger model compared to a
corresponding 1-TMLE. On one hand, as explicitly reflected in Theorem
\ref{thsecondordertmleex}, for example, it is generally true that whenever a
1-TMLE is efficient, so will be a 2-TMLE. This illustrates that 2-TMLE operates
in a safe haven wherein we expect not to hurt a 1-TMLE by performing the
additional targeting required to construct a 2-TMLE. On the other hand, we note
that 2-TMLE will be efficient in many instances in which 1-TMLE is not. As an
illustration, suppose in the setting of our motivating example that $W$ is a
univariate random variable with a sufficiently smooth density function. Suppose
also that $g_0$ is smooth enough so that an optimal univariate second-order
kernel smoother can be utilized to produce an estimate of $g_0$. In this case,
efficiency of a 1-TMLE requires that $\hat{\bar Q}$ tends to $Q_0$ at a rate
faster than $n^{-1/10}$. In contrast, the corresponding 2-TMLE built upon a
second-order canonical gradient approximated using an optimal second-order
kernel smoother will be efficient provided that $\hat{\bar Q}$ is consistent for
$\bar{Q}_0$, irrespective of the actual rate of convergence. The difference
between these requirements may not seem drastic in settings where $\bar{Q}_0$ is
sufficiently smooth since then constructing an estimator $\hat{\bar Q}$ which
satisfies both requirements is easy. This is certainly not so if $\bar{Q}_0$
fails to be smooth, in which case achieving convergence even at $n^{-1/10}$-rate
may be a challenge. This problem is exacerbated further if $W$ has several
components. For example, if $W$ is 5-dimensional, a 1-TMLE requires that
$\hat{\bar Q}$ tend to $\bar{Q}_0$ faster than $n^{-5/18}$, whereas the
corresponding 2-TMLE based on a third-order kernel-smoothed approximation
requires that $\hat{\bar Q}$ tend to $\bar{Q}_0$ faster than $n^{-1/5}$. While
the latter is achievable using an optimal second-order kernel smoother, the
former is not, and without further smoothness assumptions on $\bar{Q}_0$, a
1-TMLE will generally not be efficient.


\paragraph{Comparison with Alternative Second-Order Estimators.}

To the best of our knowledge, the only second-order estimator
preceding our proposal is discussed in \cite{robins2009quadratic}. For
a fixed bandwidth $h$, their estimator is defined as
\begin{equation}
  \hat\psi_h = \Psi(\hat P) + P_n D^{(1)}(\hat P) + \frac{1}{2} P_n^2
  D^{(2)}_h(\hat P).\label{Restim}
\end{equation}
Unlike our proposal, this estimator involved direct computation of $D^{(2)}_h$,
which in turn involves inverse weighting by an estimated multivariate density
estimate $\hat q_W(w)$. As a consequence of the curse of dimensionality these
weights may be very unstable, which may lead to a highly variable estimator in
practice. In addition, the above estimator does not always satisfy the global
constraints on the parameter space. In contrast, our proposed 2-TMLE is always
in the parameter space, since it is defined as a substitution estimator.

\subsection{Second-Order Estimator with Kernel Smoothing on the
  Missingness Score}\label{estimag}

As transpires from the developments above, even if the support of $W$ is finite
but nonetheless rich, large samples will be required to ensure that the
non-parametric estimator behaves sufficiently well. Given the sufficiency
property of the propensity score as a summary of potential confounders, it is
natural to inquire whether the use of a second-order partial gradient based on
the propensity score (see discussion in \cite{carone2014}) may allow us to
circumvent the dimensionality of $W$. Suppose that $W$ is finitely supported,
and consider the second-order expansion (\ref{soexp1}) with
\begin{align*}
  D^{(2)}(P)(o_1,o_2) &=
  \frac{2a_1\one\{g_0(w_1)=g_0(w_2)\}}{g(w_1)q_W(w_1)}\left\{1-\frac{a_2}{g(w_1)}\right\}\{y_1
  - \bar Q(w_1)\},\\
  R_3(P,P_0)&=\int\left\{1-\frac{g_0(w)q_{W,0}(w)}{g(w)q_W(w)}\right\}\left\{1-\frac{g_0(w)}{g(w)}\right\}\left\{\bar
    Q(w)-\bar Q_0(w)\right\}\dd Q_{W,0}(w).
\end{align*}
In contrast to the previous section, here $q_{W,0}(w)$ represents the density
function $\frac{d}{dx}P_0(g_0(W) \leq x)|_{x=g_0(w)}$, and $q_W(w)$ represents
$\frac{d}{dx}P(g_0(W) \leq x)|_{x=g_0(w)}$. Analogous to the multivariate case,
it is often necessary to consider a kernel function $K_h(g_0(w_1)-g_0(w_2))$
instead of the indicator $\one\{g_0(w_1)=g_0(w_2)\}$, which may not be well
supported in the data. We again denote the approximate second-order influence
function obtained with such an approximation by $D^{(2)}_h$ to emphasize the
dependence on the choice of bandwidth. Using this approximation the estimation
procedure described in the previous section may be carried out in exactly the
same fashion, but with $\hat g_h$ replaced by
\[\hat g_h(w)=\frac{\sum_{i=1}^nK_h(g_0(w) - g_0(W_i))A_i}{\sum_{i=1}^nK_h(g_0(w) - g_0(W_i))}.\]
This algorithm yields an asymptotically linear estimator of $\psi_0$ under the
assumption that $R_3(\hat P, P)= o_P(1/\sqrt{n})$, among other
regularity assumptions.

Since $g_0$ is often unknown, we must instead use an estimate $\hat g$ of $g_0$;
for example, we may take:
\[\hat g_h(w):=\frac{\sum_{i=1}^nK_h(\hat g(w) - \hat g(W_i))A_i}{\sum_{i=1}^nK_h(\hat
  g(w) - \hat g(W_i))}.\] Unfortunately, a careful analysis of the remainder
term associated with this estimator reveals that the introduction of an estimate
$\hat g$ in place of $g_0$ yields a second-order remainder term. This implies
that asymptotic efficiency of this estimator, denoted $1^\star$-TMLE, requires a
second-order term to be negligible in order to be $\sqrt{n}$-consistent. The
second-order term associated to this $1^\star$-TMLE, however, is different from
$R_2$ defined in (\ref{defR2}) and required for asymptotic linearity of the
1-TMLE. As a consequence, these estimators are expected to have different finite
sample properties. We conjecture that the $1^\star$-TMLE of this section has
improved finite sample properties over the 1-TMLE, and present a case study in
Section \ref{simula} supporting our conjecture.

\section{Simulation Studies}\label{simula}

In this section we present the results of two simulation studies,
illustrating the improvements obtained by the $1^\star$-TMLE and
2-TMLE compared to the 1-TMLE. We use covariate dimensions $d=1$ and
$d=3$ and sample sizes $n\in \{500,
1000, 2000, 10000\}$ to assess the performance of the estimators in different
scenarios. Kernel smoothers were computed using the R package ks \cite{ks}. The
bandwidth was chosen using the default method of that package
\cite{wand1994multivariate}.

\subsection{Simulation Study with $d=1$}\label{simw}

\paragraph{Simulation Setup} For each sample size $n$, we simulated 1000
datasets from the joint distribution implied by the conditional distributions
\begin{align*}
  W&\sim 6\times Beta(1/2, 1/2) - 3\\
  A|W&\sim Ber(\expit(1 + 0.7\times W))\\
  Y|A=1,W&\sim Ber(\expit(-3 +0.5\times\exp(W) + 0.5\times W)),
\end{align*}
where $Ber(\cdot)$ denotes the Bernoulli distribution, $\expit$ denotes the
inverse of the $\logit$ function, and $Beta(a,b)$ denotes the Beta distribution.

For each dataset, we fitted correctly-specified parametric models for $\bar Q_0$
and $g_0$. For a perturbation parameter $p$, we then varied the convergence rate
of $\hat{\bar Q}$ by multiplying the linear predictor by a random variable with
distribution $U(1-n^{-p}, 1)$ and subtracting a Gaussian random variable with
mean $3\times n^{-p}$ and standard deviation $n^{-p}$. Analogously, the
convergence rate of $\hat g$ was varied using a perturbation parameter $q$ by
multiplying the linear predictor by a random variable $U(1-n^{-q}, 1)$ and
subtracting a Gaussian random variable with mean $3\times n^{-q}$ and standard
deviation $n^{-q}$. We varied the values of $p$ and $q$ in a grid $\{0.01, 0.02,
0.05, 0.1, 0.2, 0.5\}^2$. This perturbation of the MLE in a correctly specified
parametric models is carried out to obtain initial estimators that have varying
consistency rates. This allows us to easily vary the convergence rate in order
to assess the performance of the estimators under such scenarios. To see how we
this procedure achieves varying consistency rates, denote the MLE of $g_0$ in
the correct parametric model by $\hat g^{\tiny \mbox{MLE}}$, and denote the
perturbed estimate by $\hat g^{\tiny \mbox{MLE}}_q$. Let $U_n$ and $V_n$ be
random variables distributed $U(1-n^{-q}, 1)$ and $N(-3n^{-q},n^{-2q})$,
respectively. Then, substituting $\hat g^{\tiny \mbox{MLE}}_q(W) = \hat g^{\tiny
  \mbox{MLE}}(W)U_n + V_n$ into $||\hat g^{\tiny \mbox{MLE}}_q - g_0||_{P_0}^2$
yields
\begin{align*}
||\hat g^{\tiny \mbox{MLE}}_q - g_0||_{P_0}^2 &\leq ||U_n(\hat g^{\tiny
  \mbox{MLE}} - g_0)||_{P_0}^2 + ||g_0(U_n-1)||_{P_0}^2 + ||V_n||_{P_0}^2\\
&= O_P(n^{-1} + n^{-2q}).
\end{align*}
Consider now different values of $q$. For example, $q=0.5$ yields the
parametric consistency rate $||\hat g^{\tiny \mbox{MLE}}_q - g_0||_{P_0}^2 =
O_P(1/n)$, whereas $q = 0$ yields an inconsistent estimator.

We computed a 1-TMLE, $1^\star$-TMLE, as well as a 2-TMLE for each initial
estimator $(\hat{\bar Q}, \hat g)$ obtained through this perturbation. We
compare the performance of the two estimators through their bias inflated by a
factor $\sqrt{n}$, relative variance compared to the nonparametric efficiency
bound, and the coverage probability of $95\%$ confidence interval assuming a
known variance. We assume the variance is known in order to isolate randomness
and bias in its estimation. The variance, bias, and coverage probabilities are
approximated through empirical means across the 1000 simulated datasets.

\paragraph{Simulation Results}

\begin{figure}[H]
  \centering
  \caption{Absolute bias of the estimators (multiplied by $\sqrt{n}$) for
    different sample sizes and convergence rates of the initial estimators of
    $\bar Q_0$ and $g_0$, when $d=1$.}
  \label{bias1}
  \includegraphics[scale=0.5]{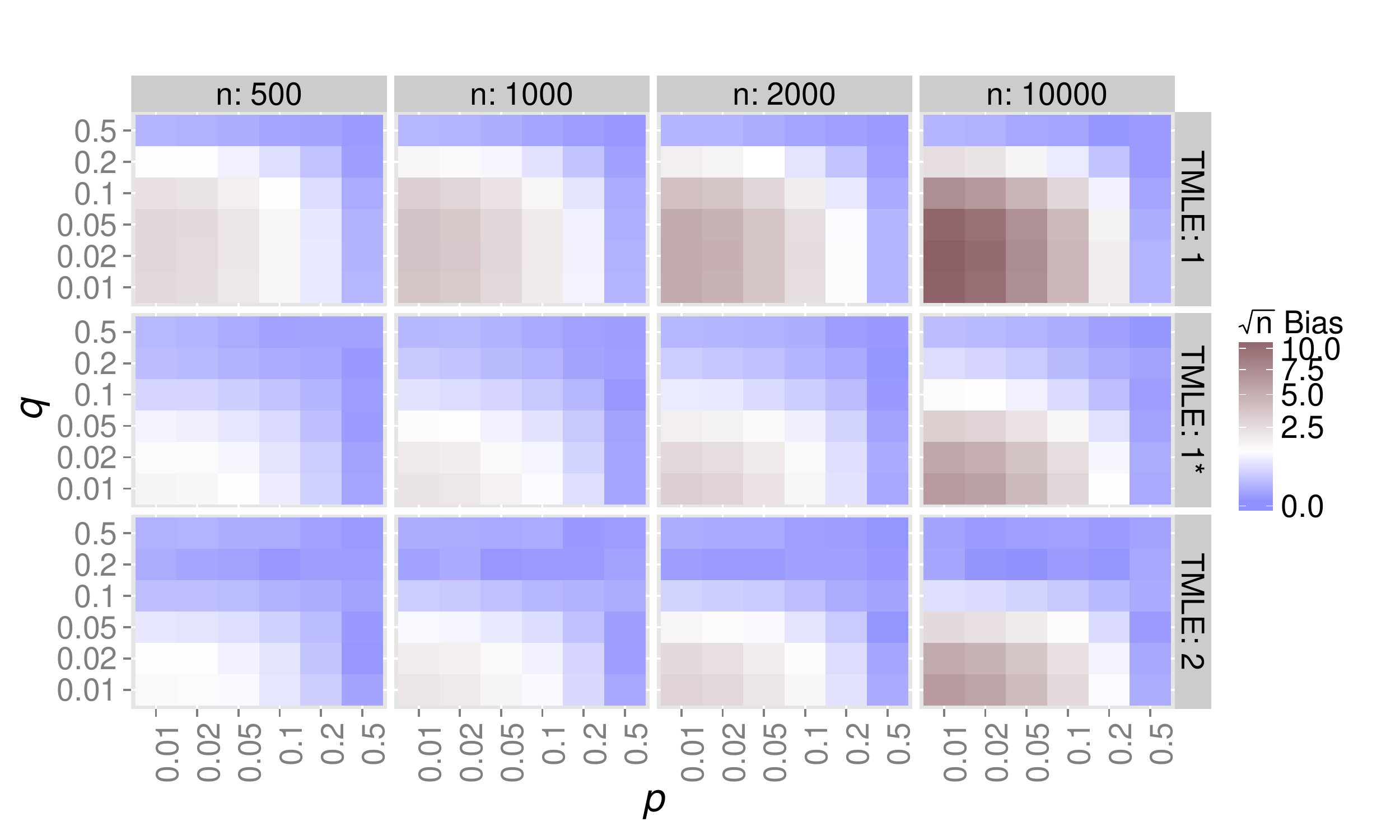}
\end{figure}

\begin{figure}[H]
  \centering
  \caption{Coverage probabilities of confidence intervals for different sample
    sizes and varying convergence rates of the initial estimators of $\bar Q_0$
    and $g_0$, when $d=1$.}
  \label{cov1}
  \includegraphics[scale=0.5]{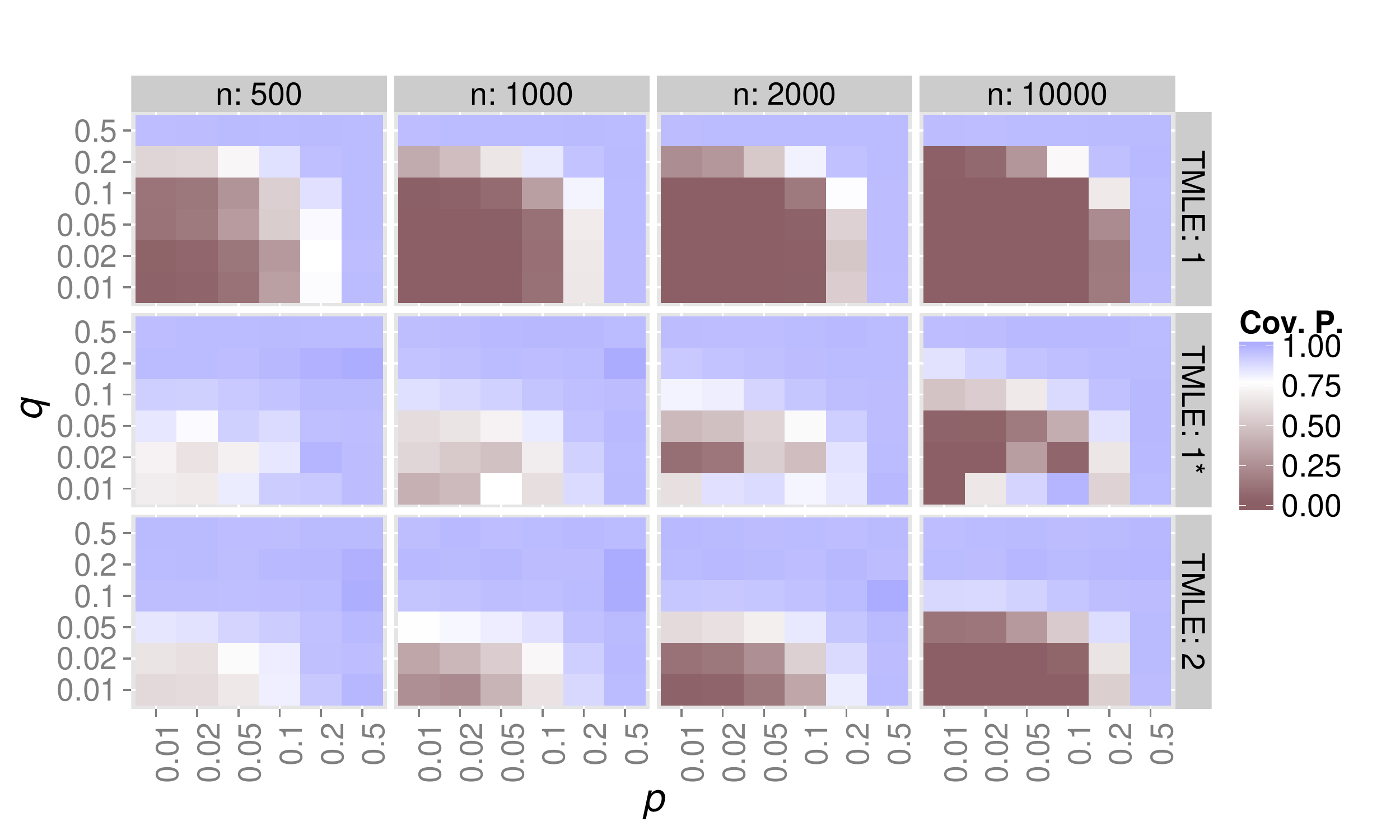}
\end{figure}

Table~\ref{stats1} shows the relative variance (rVar, defined as $n$ times the
variance divided by the efficiency bound), the absolute bias inflated by a
factor $\sqrt{n}$, as well as the coverage probability of a $95\%$ confidence
interval using the true variance of the estimators for selected values of the
perturbation parameter $(p,q)$.  Figure~\ref{bias1} shows the absolute bias of
each estimator multiplied by $\sqrt{n}$, and Figure~\ref{cov1} shows the
coverage probability of a $95\%$ confidence interval.

\setlength{\tabcolsep}{3pt}
\begin{table}[H]
  \centering
  \caption{Performance of the estimators for different sample sizes and convergence rates of the initial estimators of $\bar Q_0$ and $g_0$, when $d=1$.}
  \label{stats1}{\small
  \begin{tabular}{ccc|rrrr|rrrr|rrrr}
    \hline
  &  &       & \multicolumn{4}{|c}{1-TMLE} &    \multicolumn{4}{|c}{$1^\star$-TMLE} & \multicolumn{4}{|c}{2-TMLE}   \\\hline
  &  &       & \multicolumn{4}{|c}{$n$}    & \multicolumn{4}{|c}{$n$} & \multicolumn{4}{|c}{$n$}                         \\
  &  $p$                    & $q$   & 500                         & 1000 & 2000 & 10000 & 500  & 1000 & 2000 & 10000 & 500  & 1000 & 2000 & 10000 \\
    \hline
          &      & 0.01 & 2.43 & 3.44 & 4.86 & 10.76 & 1.31 & 1.92 & 2.84 & 5.98 & 1.19 & 1.87 & 2.66 & 5.94\\
          & 0.01 & 0.10 & 2.06 & 2.79 & 3.69 & 6.93  & 0.38 & 0.54 & 0.67 & 1.10 & 0.17 & 0.29 & 0.35 & 0.49\\
          &      & 0.50 & 0.10 & 0.11 & 0.11 & 0.10  & 0.13 & 0.12 & 0.12 & 0.15 & 0.08 & 0.07 & 0.06 & 0.03\\ \cline{3-15}
          &      & 0.01 & 1.25 & 1.65 & 2.19 & 4.15  & 0.69 & 0.95 & 1.26 & 2.45 & 0.61 & 0.91 & 1.25 & 2.41\\
  $\sqrt{n}$$|$Bias$|$    & 0.10 & 0.10 & 1.03 & 1.30 & 1.61 & 2.48  & 0.20 & 0.26 & 0.29 & 0.45 & 0.09 & 0.12 & 0.16 & 0.26\\
          &      & 0.50 & 0.04 & 0.04 & 0.03 & 0.03  & 0.03 & 0.05 & 0.07 & 0.07 & 0.06 & 0.06 & 0.02 & 0.02\\ \cline{3-15}
          &      & 0.01 & 0.11 & 0.11 & 0.10 & 0.10  & 0.03 & 0.04 & 0.04 & 0.05 & 0.03 & 0.04 & 0.05 & 0.07\\
          & 0.50 & 0.10 & 0.06 & 0.06 & 0.05 & 0.03  & 0.02 & 0.01 & 0.02 & 0.01 & 0.02 & 0.06 & 0.01 & 0.05\\
          &      & 0.50 & 0.01 & 0.00 & 0.01 & 0.01  & 0.03 & 0.01 & 0.00 & 0.00 & 0.01 & 0.01 & 0.00 & 0.03\\\hline
          &      & 0.01 & 1.42 & 1.41 & 1.45 & 1.36  & 1.74 & 1.92 & 8.37 & 3.29 & 1.75 & 1.80 & 1.70 & 1.79\\
          & 0.01 & 0.10 & 1.61 & 1.56 & 1.52 & 1.42  & 1.17 & 1.28 & 1.24 & 1.18 & 1.32 & 1.25 & 1.27 & 1.13\\
          &      & 0.50 & 1.10 & 1.11 & 1.10 & 1.09  & 1.10 & 1.14 & 1.13 & 1.12 & 1.12 & 1.12 & 1.15 & 1.16\\ \cline{3-15}
          &      & 0.01 & 1.00 & 0.98 & 0.96 & 0.94  & 1.24 & 1.16 & 1.26 & 1.33 & 1.13 & 1.14 & 1.10 & 1.05\\
  rVar    & 0.10 & 0.10 & 1.18 & 1.10 & 1.05 & 0.99  & 1.04 & 1.02 & 1.04 & 0.99 & 1.14 & 0.98 & 0.96 & 1.09\\
          &      & 0.50 & 0.97 & 1.00 & 0.97 & 0.97  & 1.00 & 1.01 & 0.96 & 0.87 & 1.09 & 0.97 & 0.97 & 1.05\\ \cline{3-15}
          &      & 0.01 & 1.03 & 1.00 & 1.04 & 1.02  & 1.04 & 1.05 & 1.00 & 0.97 & 1.03 & 1.08 & 1.05 & 0.95\\
          & 0.50 & 0.10 & 1.00 & 1.00 & 0.97 & 1.02  & 0.93 & 0.97 & 0.93 & 0.99 & 2.56 & 4.55 & 1.02 & 0.96\\
          &      & 0.50 & 0.99 & 0.98 & 0.95 & 0.99  & 0.97 & 0.96 & 0.91 & 1.00 & 0.96 & 1.00 & 1.01 & 1.02\\\hline
          &      & 0.01 & 0.02 & 0.00 & 0.00 & 0.00  & 0.52 & 0.21 & 0.60 & 0.00 & 0.56 & 0.22 & 0.02 & 0.00\\
          & 0.01 & 0.10 & 0.10 & 0.01 & 0.00 & 0.00  & 0.89 & 0.84 & 0.78 & 0.47 & 0.94 & 0.92 & 0.91 & 0.86\\
          &      & 0.50 & 0.94 & 0.94 & 0.94 & 0.94  & 0.94 & 0.94 & 0.94 & 0.94 & 0.95 & 0.94 & 0.95 & 0.94\\ \cline{3-15}
          &      & 0.01 & 0.30 & 0.09 & 0.01 & 0.00  & 0.77 & 0.59 & 0.39 & 0.00 & 0.79 & 0.61 & 0.33 & 0.00\\
  Cov. P. & 0.10 & 0.10 & 0.52 & 0.29 & 0.12 & 0.00  & 0.92 & 0.92 & 0.91 & 0.85 & 0.94 & 0.94 & 0.93 & 0.92\\
          &      & 0.50 & 0.95 & 0.95 & 0.95 & 0.95  & 0.95 & 0.95 & 0.95 & 0.95 & 0.94 & 0.95 & 0.95 & 0.94\\ \cline{3-15}
          &      & 0.01 & 0.94 & 0.94 & 0.94 & 0.94  & 0.94 & 0.95 & 0.95 & 0.95 & 0.96 & 0.95 & 0.95 & 0.94\\
          & 0.50 & 0.10 & 0.95 & 0.95 & 0.94 & 0.95  & 0.95 & 0.94 & 0.94 & 0.96 & 0.99 & 0.99 & 0.94 & 0.94\\
          &      & 0.50 & 0.94 & 0.95 & 0.95 & 0.95  & 0.95 & 0.94 & 0.95 & 0.95 & 0.95 & 0.94 & 0.95 & 0.95\\
    \hline
  \end{tabular}}
\end{table}

First, we notice that for certain slow convergence rates all the estimators have
a very large bias (e.g., $p=0.01$ and $q=0.01$, $p=0.1$ and $q=0.01$). In
contrast, for some other slow convergence rates, the absolute bias scaled by $\sqrt{n}$ of
the 1-TMLE diverges very fast in comparison to the 2-TMLE and $1^\star$-TMLE
(e.g., $p=0.1$ and $q=0.1$). The improvement in asymptotic absolute
bias of the proposed estimators comes at the price of increased variance in
certain small sample scenarios ($n\leq 2000$), such as when the outcome model
converges at a fast enough rate ($p=0.5$) but the missingness mechanism does not
($p=0.1$). In this case, the 1-TMLE has lower variance than its
competitors. This advantage of the first-order TMLE disappears asymptotically as
predicted by theory.

In terms of coverage, the improvement obtained with the $1^\star$-TMLE and the
2-TMLE is overwhelming for small values of both $p$ and $q$. As an example,
consider the case $n=2000$, $p=0.01$, $q=0.1$, in which the coverage probability
is 0 and $0.91$ for the 1-TMLE and the $1^\star$-TMLE, respectively.  This
simulation illustrates the potential for dramatic improvement obtained by using
the $1^\star$-TMLE and the 2-TMLE, which comes at the cost of over-coverage in
small sample sizes with a fast enough convergence rate ($n\leq 2000$, $p=0.5$,
$q=0.1$).

Figures~\ref{bias1} and \ref{cov1} show clearly a region of slow convergence
rates in which the proposed estimators outperform the standard first-order
TMLE. In addition, as seen in Figure~\ref{bias1}, we observe a small advantage
of the 2-TMLE over the $1^\star$-TMLE in terms of $\sqrt{n}$ bias.

\subsection{Simulations Study with $d=3$}\label{simg}

\paragraph{Simulation Setup} For each sample size $n\in \{500,
1000, 2000, 10000\}$, we simulated 1000 datasets from the joint
distribution implied by the conditional distributions
\begin{align*}
  W_1 &\sim Beta(2,2)\\
  W_2|W_1 & \sim Beta(2W_1, 2)\\
  W_3|W_1,W_2 &\sim Beta(2W_1, 2W_2)\\
  A|W&\sim Ber(\expit(1+0.12W_1 +0.1W_2+0.5W_3))\\
  Y|A=1,W&\sim Ber(\expit(-4 +0.2W_1 +0.3W_2 +0.5\exp(W_3))),
\end{align*}
where $Ber(\cdot)$ denotes the Bernoulli distribution, $\expit$ denotes the
inverse of the $\logit$ function, and $Beta(\cdot)$ denotes the beta
distribution. For each dataset, we fitted correctly-specified parametric models
for $\bar Q_0$ and $g_0$. We then varied the convergence rate of $\hat{\bar Q}$
and $\hat g$ by adding Gaussian random variables as in the previous subsection.

\begin{figure}[H]
  \centering
  \caption{Absolute bias of the estimators (multiplied by $\sqrt{n}$) for different sample sizes and varying
    convergence rates of the initial estimators of $\bar Q_0$ and $g_0$.}
  \label{bias}
  \includegraphics[scale=0.5]{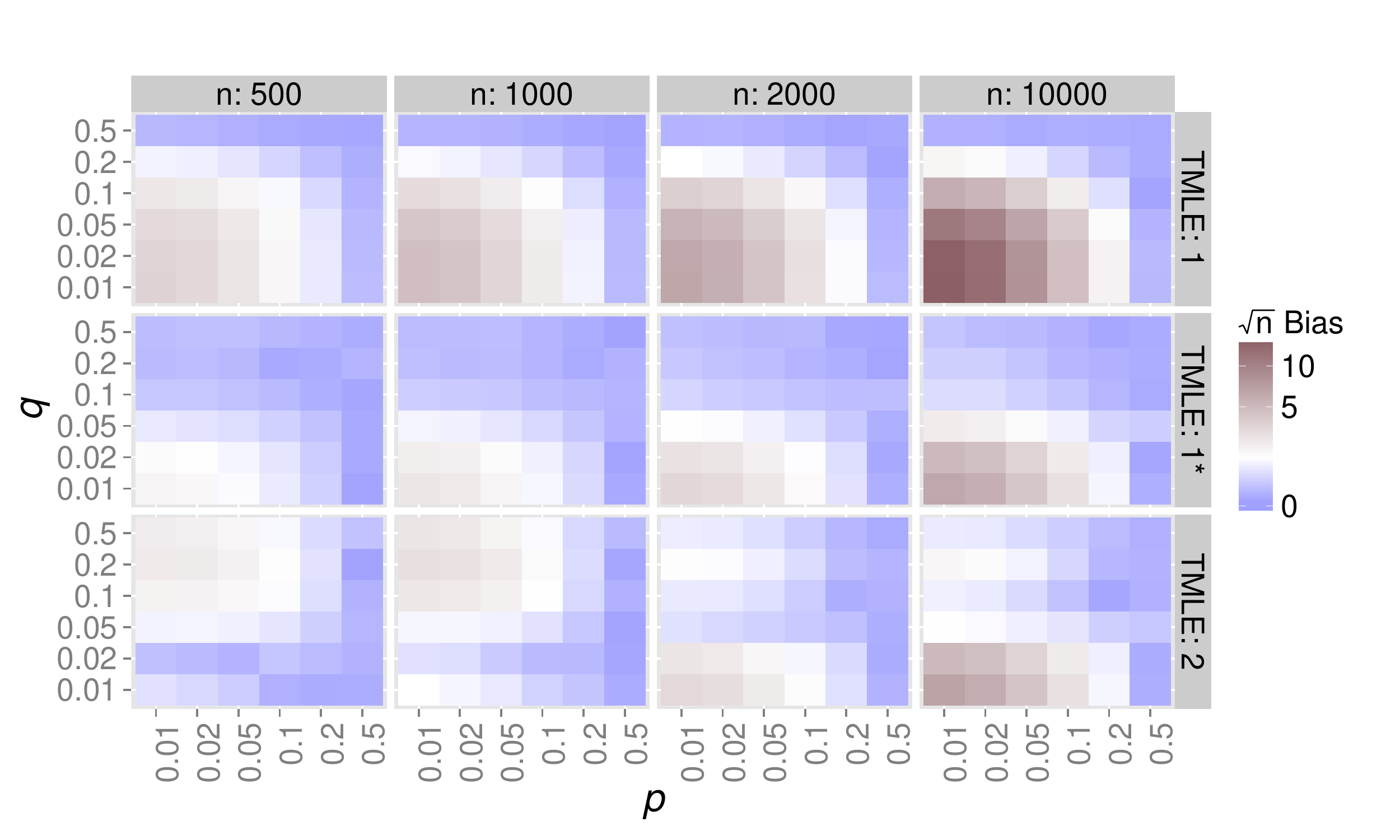}
\end{figure}

\begin{figure}[H]
  \centering
  \caption{Coverage probabilities of confidence intervals for different sample sizes and varying
    convergence rates of the initial estimators of $\bar Q_0$ and $g_0$.}
  \label{cov}
  \includegraphics[scale=0.5]{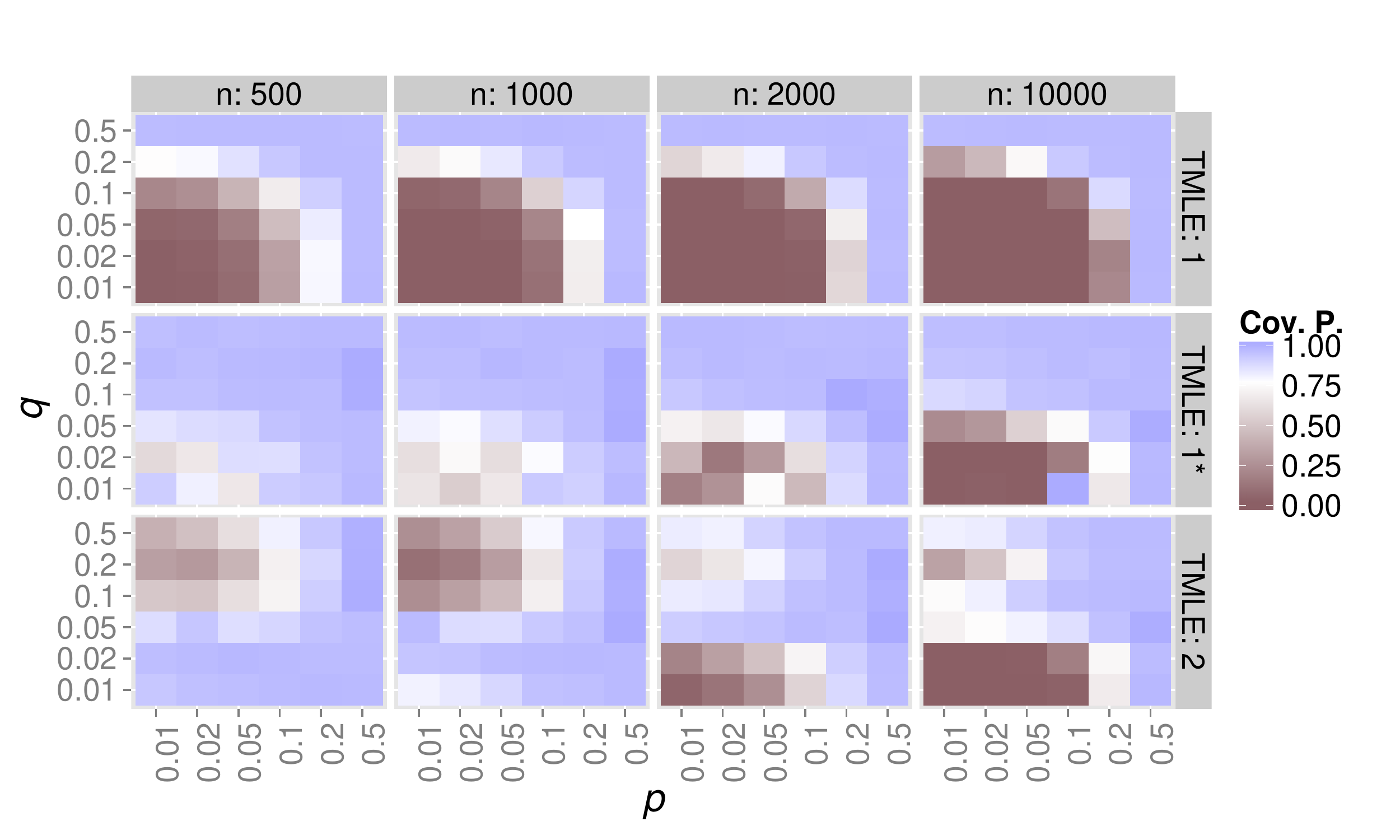}
\end{figure}

\paragraph{Simulation Results}
Table~\ref{stats2} shows the $\sqrt{n}$ absolute bias, relative variance, and
coverage probability of each estimator for selected values of the convergence
perturbation $(p,q)$. Figures~\ref{bias} and \ref{cov} show the $\sqrt{n}$
absolute bias and coverage probability of a $95\%$ confidence interval for all
values of $(p,q)$ used in the simulation.

\begin{table}[H]
  \centering
  \caption{Performance of the estimators for different sample sizes and
    convergence rates of the initial estimators of $\bar Q_0$ and $g_0$, when
    $d=3$.}
  \label{stats2}{\small
  \begin{tabular}{ccc|rrrr|rrrr|rrrr}
    \hline
  &  &       & \multicolumn{4}{|c}{1-TMLE} &    \multicolumn{4}{|c}{$1^\star$-TMLE} & \multicolumn{4}{|c}{2-TMLE}   \\\hline
  &  &       & \multicolumn{4}{|c}{$n$}    & \multicolumn{4}{|c}{$n$} & \multicolumn{4}{|c}{$n$}                         \\
  &  $p$                    & $q$   & 500                         & 1000 & 2000 & 10000 & 500  & 1000 & 2000 & 10000 & 500  & 1000 & 2000 & 10000 \\
    \hline
          &      & 0.01 & 3.02 & 4.34 & 6.14 & 13.56 & 1.39 & 1.97 & 2.77  & 5.98 & 0.47 & 1.00 & 2.69  & 6.40 \\
          & 0.01 & 0.10 & 1.93 & 2.55 & 3.27 & 5.65  & 0.18 & 0.22 & 0.32  & 0.41 & 1.47 & 1.95 & 0.61  & 0.73 \\
          &      & 0.50 & 0.07 & 0.05 & 0.04 & 0.03  & 0.09 & 0.09 & 0.12  & 0.15 & 1.71 & 2.05 & 0.67  & 0.63 \\ \cline{3-15}
          &      & 0.01 & 1.33 & 1.77 & 2.31 & 4.26  & 0.63 & 0.84 & 1.17  & 2.22 & 0.03 & 0.28 & 1.08  & 2.28 \\
  $\sqrt{n}$$|$Bias$|$    & 0.10 & 0.10 & 0.87 & 1.05 & 1.25 & 1.70  & 0.08 & 0.11 & 0.14  & 0.17 & 0.96 & 1.02 & 0.22  & 0.12 \\
          &      & 0.50 & 0.01 & 0.02 & 0.01 & 0.03  & 0.07 & 0.05 & 0.07  & 0.04 & 0.87 & 0.90 & 0.21  & 0.23 \\ \cline{3-15}
          &      & 0.01 & 0.09 & 0.08 & 0.08 & 0.07  & 0.00 & 0.01 & 0.03  & 0.03 & 0.02 & 0.02 & 0.04  & 0.02 \\
          & 0.50 & 0.10 & 0.04 & 0.03 & 0.03 & 0.00  & 0.00 & 0.07 & 0.19  & 0.01 & 0.02 & 0.09 & 0.11  & 0.03 \\
          &      & 0.50 & 0.01 & 0.00 & 0.01 & 0.01  & 0.02 & 0.00 & 0.01  & 0.01 & 0.15 & 0.08 & 0.01  & 0.03 \\\hline
          &      & 0.01 & 1.60 & 1.59 & 1.57 & 1.60  & 3.22 & 2.13 & 2.15  & 2.11 & 2.58 & 2.87 & 1.89  & 2.09 \\
          & 0.01 & 0.10 & 1.73 & 1.72 & 1.56 & 1.46  & 1.05 & 1.09 & 1.02  & 0.99 & 2.14 & 1.97 & 1.27  & 1.17 \\
          &      & 0.50 & 1.06 & 1.08 & 1.03 & 1.01  & 1.10 & 1.08 & 0.97  & 1.07 & 2.15 & 2.10 & 1.31  & 1.17 \\ \cline{3-15}
          &      & 0.01 & 1.10 & 1.07 & 1.06 & 1.04  & 1.13 & 1.26 & 1.14  & 1.25 & 1.71 & 1.58 & 1.31  & 1.19 \\
  rVar    & 0.10 & 0.10 & 1.21 & 1.20 & 1.09 & 1.04  & 0.95 & 0.99 & 0.99  & 0.92 & 1.75 & 1.78 & 1.17  & 1.06 \\
          &      & 0.50 & 0.97 & 0.98 & 1.02 & 0.98  & 1.03 & 1.01 & 1.01  & 1.01 & 2.17 & 1.96 & 1.19  & 1.04 \\ \cline{3-15}
          &      & 0.01 & 1.02 & 1.04 & 1.01 & 1.01  & 1.04 & 1.00 & 1.05  & 1.03 & 1.02 & 1.02 & 0.99  & 0.99 \\
          & 0.50 & 0.10 & 1.04 & 0.95 & 0.97 & 0.97  & 1.14 & 6.10 & 17.17 & 1.02 & 3.58 & 9.19 & 16.22 & 0.84 \\
          &      & 0.50 & 0.99 & 0.98 & 1.01 & 0.95  & 1.10 & 0.93 & 1.00  & 0.93 & 1.37 & 1.23 & 1.10  & 0.97 \\\hline
          &      & 0.01 & 0.01 & 0.00 & 0.00 & 0.00  & 0.76 & 0.23 & 0.03  & 0.00 & 0.91 & 0.78 & 0.03  & 0.00 \\
          & 0.01 & 0.10 & 0.19 & 0.03 & 0.00 & 0.00  & 0.93 & 0.92 & 0.90  & 0.86 & 0.49 & 0.22 & 0.80  & 0.73 \\
          &      & 0.50 & 0.94 & 0.95 & 0.94 & 0.94  & 0.93 & 0.95 & 0.95  & 0.94 & 0.38 & 0.22 & 0.80  & 0.79 \\ \cline{3-15}
          &      & 0.01 & 0.30 & 0.08 & 0.01 & 0.00  & 0.79 & 0.69 & 0.42  & 0.03 & 0.95 & 0.93 & 0.54  & 0.02 \\
  Cov. P. & 0.10 & 0.10 & 0.66 & 0.53 & 0.35 & 0.10  & 0.95 & 0.94 & 0.94  & 0.93 & 0.70 & 0.67 & 0.93  & 0.94 \\
          &      & 0.50 & 0.95 & 0.95 & 0.95 & 0.94  & 0.94 & 0.95 & 0.94  & 0.94 & 0.79 & 0.77 & 0.92  & 0.93 \\ \cline{3-15}
          &      & 0.01 & 0.95 & 0.95 & 0.95 & 0.95  & 0.95 & 0.95 & 0.96  & 0.96 & 0.95 & 0.95 & 0.94  & 0.96 \\
          & 0.50 & 0.10 & 0.94 & 0.94 & 0.95 & 0.95  & 0.95 & 0.99 & 0.98  & 0.95 & 0.99 & 0.98 & 0.98  & 0.95 \\
          &      & 0.50 & 0.95 & 0.95 & 0.95 & 0.95  & 0.94 & 0.95 & 0.95  & 0.96 & 0.93 & 0.95 & 0.95  & 0.95 \\
    \hline
  \end{tabular}}
\end{table}

The remarks of the previous section regarding the trade-offs between variance
and bias in different regions of the convergence rates also hold for this
simulation. The main difference observed here is that the 2-TMLE has poorer
performance in terms of $\sqrt{n}$ bias than the 1-TMLE and the $1^\star$-TMLE
for small samples when one of the models converges at a fast enough rate
($p=0.5$ or $q=0.5$). This problem disappears somewhat as $n$ increases, but
it highlights the point that the 2-TMLE should be used with caution in small
samples.

In this simulation we do not see any practical advantage of the 2-TMLE over the
$1^\star$-TMLE. In fact, the $1^\star$-TMLE performs better than the 2-TMLE for
small samples, and outperforms the 1-TMLE in all sample sizes, with the caveat
of increased variance in certain scenarios as discussed in the previous
section.

\section{Data Illustration}\label{aplica}

In order to illustrate the method presented, we make use of the dataset
\texttt{lindner} available in the R package \texttt{PSAgraphics}. The
dataset contains data on 996 patients treated at the Lindner
Center, Christ Hospital, Cincinnati in 1997, and were originally
analyzed in \cite{bertrand2002management}. All patients received a Percutaneous Coronary
Intervention (PCI). One of the primary goals of the original study was
to assess whether administration of Abciximab, an anticoagulant,
during PCI improves short and long term health outcomes of
patients undergoing PCI. We reanalyze the \texttt{lindner} dataset
focusing on the cardiac related costs incurred within 6
months of patient’s initial PCI as an outcome. The covariates measured
are: indicator of coronary stent deployment during the PCI,
height, sex, diabetes status, prior acute myocardial infarction, left
ejection fraction, and number of vessels involved in the PCI.

As noted by several authors \cite[e.g.][]{Rosenbaum&Rubin83,
  bang2005doubly, mohan2013missing}, causal inference problems may be
tackled using methods for missing data. Let $T$ denote
an indicator of having received Abciximab. Adopting the potential
outcomes framework, consider the potential outcomes $Y_t$, $t\in\{0,1\}$,
given by the outcomes that would have been observed in a hypothetical world if,
contrary to fact, $P(T=t)=1$. The consistency assumption
states that $A=t$ implies that $Y_t=Y$, where $Y$ is the observed
outcome. Thus, $E(Y_t)$ may be estimated using methods for missing
outcomes, where $Y_t$ is observed only when $T= t$. In particular,
estimation of $E(Y_1)$ and $E(Y_0)$ is carried out using the methods
described in the previous sections with $A=T$ and $A=1-T$,
respectively. Our parameter of interest is the average treatment
effect $E(Y_1)-E(Y_0)$.

Since the outcome is continuous, we first used the transformation
$(y-\min(y))/(\max(y)-\min(y))$ to map it to the interval $[0,1]$. We
then used the approach outlined in \cite{gruber2010targeted} to
construct the 1-TMLE and the $1^\star$-TMLE. The distribution of both
estimators was estimated with the bootstrap as discussed in Section~4.2 of
\cite{carone2014}, which involves bootstrapping the second-order expansion of
the estimator. This bootstrapped distribution is expected to capture the
second-order behavior of the estimators and therefore possibly results in finite
sample gains for the $1^\star$-TMLE.

The mean of the outcome conditional on covariates was estimated separately for
the two treatment groups. Both the outcome regression and the treatment
mechanism where estimated using a model stacking technique called Super Learning
\cite{vanderLaan&Polley&Hubbard07}. Super Learning takes a collection of
candidate estimators and combines them in a weighted average, where the weights
are chosen to minimize the cross-validated prediction error of the final
predictor, measured in terms of the $L^2$ loss function. The collection of
algorithms used is described in Table~\ref{algs}.  Table~\ref{sl} shows the
cross-validated risks of the algorithms as well as their weights in the final
predictor of $\bar Q_0$ and
$g_0$. 

\begin{table}[H]
  \centering
  \caption{Prediction algorithms used to estimate $\bar Q_0$ and $g_0$}\label{algs}
  \begin{tabular}{r p{10cm}}\hline
    Algorithm & Description\\\hline
    GLM & Generalized linear model. The $\logit$ link was used for $g_0$
    and the identity for $Q_0$.\\\hline
    BayesGLM & Bayesian GLM. Weakly informative priors were used as implemented
    by default in the function \texttt{bayesglm} of the \texttt{arm}
    package in R.\\\hline
    GAM & Generalized additive model as implemented in the R package
    \texttt{gam}.\\\hline
    PolyMARS & Multivariate adaptive polynomial spline regression
    implemented in the R package \texttt{polspline}.\\\hline
    Earth & Multivariate adaptive regression splines implemented in the R
    package \texttt{earth}.\\\hline
  \end{tabular}
\end{table}

\begin{table}[H]
  \centering
  \caption{Cross-validated risk and weight of each algorithm in the
    Super Learner for estimation of $\bar Q_0$ and $g_0$.}\label{sl}
  \begin{tabular}{l|rr|rr|rr}    \hline
    & \multicolumn{2}{c}{$\bar Q_0$ Treated} & \multicolumn{2}{|c}{$\bar Q_0$ Untreated} & \multicolumn{2}{|c}{$g_0$}           \\\hline
    Algorithm & CV Risk                               & Weight                                    & CV Risk & Weight  & CV Risk & Weight  \\\hline
    GLM                 & 0.00275                               & 0.00000                                   & 0.00684 & 0.00000 & 0.19506 & 0.00000 \\
    BayesGLM            & 0.00275                               & 0.00000                                   & 0.00684 & 0.00000 & 0.19502 & 0.13993 \\
    GAM                 & 0.00274                               & 0.65699                                   & 0.00679 & 0.57261 & 0.19495 & 0.00000 \\
    PolyMARS            & 0.00280                               & 0.15156                                   & 0.00709 & 0.21333 & 0.18905 & 0.62503 \\
    Earth               & 0.00281                               & 0.19145                                   & 0.00688 & 0.21405 & 0.19332 & 0.23504 \\\hline
  \end{tabular}
\end{table}

For bandwidth selection, we use a loss function that targets directly the
first-order expansion of the parameter of interest, which is equivalent to the
first step of the collaborative TMLE (C-TMLE) presented in
\cite{vanderLaan&Gruber09}. This approximation of the C-TMLE is computationally
more tractable and is justified theoretically as argued below.

Following \cite{gruber2011ctmle}, let $s\in\{1,\ldots,S\}$ index a random sample
split into a validation sample $V(s)$ and a training sample $T(s)$. The
cross-validation bandwidth selector is defined as
\[\hat h:=\argmin_h\{cvRSS(h) + cvVar(h) + n\times [cvBias(h)]^2\},\]
where
\begin{align*}
  cvRSS(h) &:= \sum_{s=1}^S\sum_{i \in V(s)}\{Y_i - \hat{\bar
    Q}^\star_{h,s}(W_i)\}^2,\\
  cvVar(h)&:= \sum_{s=1}^S\sum_{i \in V(s)}\left[\frac{A_i}{\hat
      g_s(W_i)}\{Y_i - \hat{\bar Q}^\star_{h,s}(W_i)\} + \hat{\bar
      Q}^\star_{h,s}(W_i) - \hat\psi_{h,s}\right]^2,\text{ and }\\
  cvBias(h)&:= \frac{1}{S}\sum_{s=1}^S\left(\hat\psi_{h,s} - \hat\psi_{h}\right)
\end{align*}
are the cross-validated residual sum of squares (RSS), cross-validated variance
estimate, and cross-validated bias estimate, respectively. The key idea is to
select the bandwidth $h$ that makes $\hat H^{(2)}_h$ most predictive of $Y$,
while adding an asymptotically negligible penalty term for increases in bias and
variance in estimation of $\psi_0$. Here, $\hat{\bar Q}^\star_{h,s}$,
$\hat\psi_{h,s}$, and $\hat g_s$ are the result of applying the estimation
algorithms described in Section \ref{soet} using only data in the training
sample $T(s)$. 

This loss function is the result of adding a mean squared error (MSE) term
$cvVar(h) + n\times [cvBias(h)]^2$ to the usual RSS loss function
used in regression problems. Since he MSE contribution to the loss function
is asymptotically negligible compared to the RSS, this loss function
yields a valid loss function for the parameter $\bar Q_0$. Intuitively, the
cross-validated MSE term serves the purpose of penalizing bandwidths that are
solely targeted to estimation of $\bar Q_0$ but perform poorly for
$\psi_0$. This bandwidth selection algorithm as well as the estimator is implemented in
the R code provided in the supplementary materials.

\paragraph{Results} The unadjusted dollar difference in the outcome between the
two groups is equal to US\$$1512$. The 1-TMLE and the $1^\star$-TMLE give an
adjusted difference of US\$$765$ and US\$$561$; with $95\%$ confidence intervals
$(-667, 2732)$ and $(-1212, 2174)$, respectively. The bootstrap standard errors
of the two estimators are $803$ and $826$, respectively. The larger variance of
the $1^\star$-TMLE may be an consequence of our conjectured property that the
$1^\star$-TMLE has a better finite sample bias-variance trade-off. In this
illustration the use of an estimator with improved asymptotic properties
considerably changes the point estimate and confidence intervals.

\section{Discussion}
We proposed a second-order estimator of the mean of an outcome missing at
random, and present a theorem showing the conditions under which it is expected
to be asymptotically efficient. Our main accomplishment is to show that the
second-order TMLE achieves efficiency under slower convergence rates of the
initial estimators than those required for efficiency of first-order
estimators. The conditions for effficiency of our proposed second-order
procedure include the convergence of a kernel bandwidth estimator at rates that
are not allowed to be too fast or too slow. The construction of algorithms that
achieve the required rates remains an open question.

In addition to the second-order estimator, we presented a novel first-order
estimator whose construction is inspired by a second-order expansion of the
parameter functional. We showed dramatic improvements in bias and coverage
probability of this estimator compared to a first-order competitor in
simulations. We conjecture that gains of this kind are expected to hold in
general for finite samples, but a formal study of the remainder terms of
both estimators remains to be done.

\section{Acknowledgments}
Marco Carone was supported in part by NIH grant 2R01AI074345 and by an endowment
generously provided by Genentech. Mark J. van der Laan was supported NIH grant
R01AI07434506.

\bibliographystyle{plainnat}
\bibliography{tmle}

\begin{thebibliography}{18}
\providecommand{\natexlab}[1]{#1}
\providecommand{\url}[1]{\texttt{#1}}
\expandafter\ifx\csname urlstyle\endcsname\relax
  \providecommand{\doi}[1]{doi: #1}\else
  \providecommand{\doi}{doi: \begingroup \urlstyle{rm}\Url}\fi

\bibitem[Bang and Robins(2005)]{bang2005doubly}
Heejung Bang and James~M Robins.
\newblock Doubly robust estimation in missing data and causal inference models.
\newblock \emph{Biometrics}, 61\penalty0 (4):\penalty0 962--973, 2005.

\bibitem[Bertrand et~al.(2002)Bertrand, Simoons, Fox, Wallentin, Hamm, Feyter,
  Specchia, Ruzyllo, and McFadden]{bertrand2002management}
Micheal~E Bertrand, Maarten~L Simoons, Keith~AA Fox, Lars~C Wallentin,
  Christian~W Hamm, PJ~Feyter, G~Specchia, Witold Ruzyllo, and EP~McFadden.
\newblock Management of acute coronary syndromes in patients presenting without
  persistent st-segment elevation.
\newblock \emph{European heart journal}, 2002.

\bibitem[Carone et~al.(2014)Carone, D{\'\i}az, and van~der Laan]{carone2014}
Marco Carone, Iv{\'a}n D{\'\i}az, and Mark~J van~der Laan.
\newblock Higher-order targeted minimum loss-based estimation.
\newblock 2014.

\bibitem[Duong(2015)]{ks}
Tarn Duong.
\newblock \emph{ks: Kernel Smoothing}, 2015.
\newblock URL \url{http://CRAN.R-project.org/package=ks}.
\newblock R package version 1.9.4.

\bibitem[Gruber and van~der Laan(2011)]{gruber2011ctmle}
Susan Gruber and Mark van~der Laan.
\newblock C-tmle of an additive point treatment effect.
\newblock In Mark van~der Laan and Sherri Rose, editors, \emph{Targeted
  Learning: Causal Inference for Observational and Experimental Data}.
  Springer, 2011.

\bibitem[Gruber and van~der Laan(2010)]{gruber2010targeted}
Susan Gruber and Mark~J van~der Laan.
\newblock A targeted maximum likelihood estimator of a causal effect on a
  bounded continuous outcome.
\newblock \emph{The International Journal of Biostatistics}, 6\penalty0 (1),
  2010.

\bibitem[Mohan et~al.(2013)Mohan, Pearl, and Tian]{mohan2013missing}
Karthika Mohan, Judea Pearl, and Jin Tian.
\newblock Missing data as a causal inference problem.
\newblock In \emph{Proceedings of the Neural Information Processing Systems
  Conference (NIPS)}. Citeseer, 2013.

\bibitem[Robins et~al.(2009{\natexlab{a}})Robins, Li, Tchetgen, and van~der
  Vaart]{robins2009quadratic}
James Robins, Lingling Li, Eric Tchetgen, and Aad~W van~der Vaart.
\newblock Quadratic semiparametric von mises calculus.
\newblock \emph{Metrika}, 69\penalty0 (2-3):\penalty0 227--247,
  2009{\natexlab{a}}.

\bibitem[Robins et~al.(2009{\natexlab{b}})Robins, Tchetgen, Li, van~der Vaart,
  et~al.]{robins2009semiparametric}
James Robins, Eric~Tchetgen Tchetgen, Lingling Li, Aad van~der Vaart, et~al.
\newblock Semiparametric minimax rates.
\newblock \emph{Electronic Journal of Statistics}, 3:\penalty0 1305--1321,
  2009{\natexlab{b}}.

\bibitem[Rubin(1983)]{Rosenbaum&Rubin83}
P.R. Rosenbaum \&~D.B. Rubin.
\newblock The central role of the propensity score in observational studies for
  causal effects.
\newblock \emph{Biometrika}, 70:\penalty0 41--55, 1983.

\bibitem[Starmans(2011)]{starmans2011models}
Richard~JCM Starmans.
\newblock Models, inference, and truth: probabilistic reasoning in the
  information era.
\newblock In Mark van~der Laan and Sherri Rose, editors, \emph{Targeted
  Learning: Causal Inference for Observational and Experimental Data}.
  Springer, 2011.

\bibitem[Tan(2014)]{tan2014second}
Zhiqiang Tan.
\newblock Second-order asymptotic theory for calibration estimators in sampling
  and missing-data problems.
\newblock \emph{Journal of Multivariate Analysis}, 131:\penalty0 240--253,
  2014.

\bibitem[van~der Laan \& E. Polley \&
  A.~Hubbard(2007)]{vanderLaan&Polley&Hubbard07}
M.J. van~der Laan \& E. Polley \& A.~Hubbard.
\newblock Super learner.
\newblock \emph{Statistical Applications in Genetics \& Molecular Biology},
  6\penalty0 (25), 2007.
\newblock ISSN 1.

\bibitem[van~der Laan and Rose(2011)]{vanderLaanRose11}
M.J. van~der Laan and S.~Rose.
\newblock \emph{Targeted Learning: Causal Inference for Observational and
  Experimental Data}.
\newblock Springer, New York, 2011.

\bibitem[van~der Laan and Rubin(2006)]{vanderLaan&Rubin06}
M.J. van~der Laan and D.~Rubin.
\newblock Targeted maximum likelihood learning.
\newblock \emph{The International Journal of Biostatistics}, 2\penalty0 (1),
  2006.

\bibitem[van~der Laan \& S.~Gruber(2009)]{vanderLaan&Gruber09}
M.J. van~der Laan \& S.~Gruber.
\newblock Collaborative double robust penalized targeted maximum likelihood
  estimation.
\newblock \emph{The International Journal of Biostatistics}, 2009.

\bibitem[van~der Vaart(1998)]{vanderVaart98}
A.~W. van~der Vaart.
\newblock \emph{Asymptotic Statistics}.
\newblock Cambridge University Press, 1998.

\bibitem[Wand and Jones(1994)]{wand1994multivariate}
MP~Wand and MC~Jones.
\newblock Multivariate plug-in bandwidth selection.
\newblock \emph{Computational Statistics}, 9\penalty0 (2):\penalty0 97--116,
  1994.

\end{thebibliography}


\end{document}